\documentclass[final,10pt,journal,letterpaper]{IEEEtran}

\IEEEoverridecommandlockouts

\usepackage{amsmath}
\usepackage{amsfonts}
\usepackage{amssymb}
\usepackage[compress]{cite}
\usepackage{graphicx}
\usepackage{algorithmicx}
\usepackage{algcompatible}
\usepackage{algorithm}
\usepackage{url}
\usepackage{color}
\usepackage{mathrsfs}
\usepackage{todonotes}
\usepackage{bbm}
\usepackage{url}
\usepackage{bm}
\usepackage{subcaption}

\newtheorem{proposition}{Proposition}
\newtheorem{assumption}{Assumption}

\newtheorem{remark}{Remark}
\newtheorem{corollary}{Corollary}
\newtheorem{definition}{Definition}

\newcommand{\QED}{}

\DeclareMathOperator*{\argmin}{arg\,min}

\DeclareMathOperator{\argmax}{arg\,max}

\newcommand{\citep}{\cite}

\begin{document}

\author{Farhad~Farokhi,~Iman~Shames, and~Michael~Cantoni
\thanks{The authors are with the Department of Electrical and Electronic Engineering, University of Melbourne, Parkville, Victoria 3010, Australia. Emails:\{ffarokhi,ishames,cantoni\}@unimelb.edu.au} }

\title{Budget-Constrained Contract Design for Effort-Averse Sensors \\ in Averaging Based
Estimation\thanks{The work was supported by a McKenzie Fellowship and the Australian Research Council (LP130100605).}}

\maketitle

\begin{abstract}  
Consider a group of effort-averse, or lazy, sensors that seek to minimize the effort invested to collect measurements of a variable. Increasing the effort invested by the sensors improves the quality of the measurements provided to the central planner but this incurs increased costs to the sensors. The central planner, which processes the sensor measurements, employs an averaging estimator. It also determines contracts for rewarding sensors based on the measurements obtained. The problem of designing a contract that yields an estimation-error based quality-of-service level in return for the reward extended to sensors is investigated in this paper. To this end, a game is formulated between the central planner and the sensors.
 Conditions for the existence and uniqueness of an equilibrium are identified. The equilibrium is constructed explicitly and its properties in response to a reward based contract are studied. It turns out that the central planner, while not being able to directly measure the effort invested by the sensors, can enhance the estimation quality by rewarding each sensor based on the distance of its measurements from the output of the averaging estimator. Ultimately, optimal contracts are designed from the perspective of the budget required for achieving a specified level of estimation error.  
\end{abstract}


\section{Introduction}
The number of networked platforms, from smart phones to wearable gadgets, has increased at a
staggering rate in the last few years. Capitalizing on the ubiquity of networked devices in various circumstances gives rise to issues that relate to the ownership of and authority over such devices. For example, crowd-sensing applications involve registered participants (or agents) providing measurements of a quantity, such as the traffic congestion in a locale or the quality of a service received from a provider, which are then combined to construct accurate estimates of this variable. This, however, opens doors to selfish behaviours when collecting the measurements due to 
resource constraints, privacy concerns, conflicts of interest, or security issues. Many studies have speculated on the use of monetary rewards to provide appropriate incentives to the participants, e.g.,~\citep{7055221}. However, there are several concerns that need to be addressed before adopting such schemes. For instance, if the rewards are not tied to the quality of the provided measurements, effort-averse\footnote{The term effort averse is borrowed from the economics literature. For instance, in insurance industry, it is observed that ``flat wage leaves the agents with no incentive to avoid behaviours ... that increase the risks'' since such avoidance requires investing effort (which is clearly not rewarded under flat wages)~\citep[p.\,352]{menard2008handbook}. Therefore, several studies have been devoted to designing optimal contracts that reward efforts either directly (based on the time and energy spent) or indirectly (based on the outcome)~\citep{christensen2005economics,menard2008handbook}.} (or lazy) agents do not have any incentive for providing high-quality measurements in view of the requisite high investment of effort (e.g. time, energy, loss of opportunity) to collect such measurements. Moreover, providing rewards that promote the sensors towards the investment of more effort is not straightforward since the quality of measurements, or equivalently the effort that is invested by the sensors, may not be fixed or known ahead of time.

Consider a group of effort-averse sensors that provide measurements of a variable of interest to a central planner. The central planner implements a simple averaging estimator and rewards the sensors to improve the quality of the collected measurements (equivalently to increases the amount invested efforts by the sensors) according to a contract that is fixed before any measurements are taken. Consideration of the simple averaging estimator is motivated by the central planner's lack of direct control over and general ignorance of the effort that individual sensors invest, whereby optimality of averaging weights cannot be reasonably defined (as in the case of the least mean square estimator). It also greatly simplifies the mathematical derivations and makes the scheme robust, as discussed below. Rewards are assigned according to contracts that entitle sensors to payments of a specific from for measurements. This paper investigates the problem of designing such contracts under budget constraints. In response to the contracts, sensors determine the level of effort to invest in order to strike a balance between the corresponding cost and the expected reward to be gained. The described interaction between the sensors and the central planner is modelled by a game. Conditions are provided to guarantee the existence and uniqueness of an equilibrium of the corresponding game. The equilibrium is constructed explicitly and its properties are investigated as a function of the parameters of the proposed contracts. Interestingly, the efforts invested by the sensors, at the captured equilibrium, are dominant strategies. That is, even if some sensors are faulty or mistakenly setting their efforts, it is in the best interest of each effort-averse sensor to expend the equilibrium effort. This property makes the presented framework robust to faults and cyber-security threats (as, in these cases, some sensors might not follow the strategies that are in their best interests). Finally, two fundamental properties for averaging estimators are proved. The first property provides the minimum required budget for achieving a specified level of the estimation error, while the second one bounds the quality of the estimates for a given budget. These fundamental properties are proved over the set of all individually rational contracts for which the \textit{ex~ante} expected return of each sensor is non-negative (as, otherwise, the sensors opt out of the sensing scheme). These fundamental properties are used to design optimal contracts under budget constraints. To provide the context for the novelty of these results, related literature on incentivizing sensors in participatory-sensing schemes is reviewed below. A brief description of the paper structure follows.

The problem of incentivizing strategic participants in crowd-sensing problems has been studied extensively in the past; see~\citep{restuccia2015incentive} for a review of the results. In~\citep{lee2010sell}, the sensors are assumed to sell homogeneous (interchangeable) sensing data to a central planner with a minimal
required level of participation. The sensors submit bids and the ones with the lowest bids are recruited and compensated. Since sensors may lose  interest in future participation if they consistently lose in the bidding process, the sensors that are not recruited in a given round are incentivized by using a virtual currency that is payable, in the future, on reduction of their submitted bids. The framework is extended in~\citep{jaimes2012location} to accommodate budget constraints. A crowd-sensing scheme for buying and selling parking information is introduced in~\citep{hoh2012trucentive}, where the rewards are designed to bring good information and deter malicious behaviour by rational individuals (e.g., adding fake data for personal gain). In~\citep{duan2014motivating}, gathering data when the sensors incur a constant cost for their reports is considered. The central planner is assumed to need more than a certain number of participants to obtain a positive revenue for itself. At the beginning of the process, the central planner announces the total reward (to be divided between the sensors) and the minimum number of sensors required. The sensors then decide to participate if the reward divided by the number of participants in that round (if the threshold is passed) is larger than their costs.  Mobile sensors moving in and out of desired sensing area are considered in~\citep{gao2015providing}. At the beginning of each time slot, a subset of the sensors are assigned to a task. The sensing range is assumed to be limited and sensing is considered to be costly. It is also assumed that sensors suffer indirect costs  when not performing any tasks. Therefore, sensors may drop out if they are not rewarded enough in the long run. The central planner aims at maximizing its return while making sure there are not many drop outs.  The central planner in~\citep{6848053} is assumed to seek sensors for completing tasks inside a region of interest. The sensors have different arrival and departure times as well as privately-known costs for each task. They bid for a group of tasks and the planner selects a group of them to perform the tasks under a budget constraint. A common feature of the schemes described above is the aim of retaining a useful set of participants while minimising expenditure of the central planner. By contrast, the work presented here focuses on influencing the underlying effort invested by the participanting sensors in order to improve the quality of estimates determined by averaging. 

This paper is closer, in essence, to the studies in~\citep{yang2012crowdsourcing,gao2015providing,6848053,han2016truthful,luo2014profit, koutsopoulos2013optimal}. In~\citep{yang2012crowdsourcing}, the sensors are assumed to incur a cost that is proportional to the time they dedicate to a sensing task. The central planner then rewards the participants by an amount that they divide among themselves according to the time spent by each to provide measurements. The objective of the central planner is to maximize a measure of the time invested by the sensors. In~\citep{han2016truthful}, the central planner wants to assign a set of tasks with different valuations among sensors under a total budget constraint. Each sensor can only perform one task at a time (but it can move between multiple tasks) and has a sensing cost per unit of time, which is not available to the central planner. The sensors also have privately-known ranges of time during which they can perform the sensing tasks. The sensors first bid prices (costs to the central planner) for their service, start times, and end times. Then, the planner announces a plan (time period allocated to tasks for each sensor) and compensate them accordingly. The goal is to maximize the total time that the sensors dedicate to the tasks. 
In~\citep{luo2014profit}, the sensors are assumed to have privately-known types (such as the marginal cost of participation or a measure of their skills) and invest efforts for gathering information. The cost of sensing is assumed to depend on the type and the effort in a specific form. The central planner only compensates the sensor that invests the highest effort (proportional to its invested effort) and the sensors maximize the expected reward minus the cost. The goal of the central planner is to maximize the total invested effort minus the reward. In~\citep{koutsopoulos2013optimal}, sensors are considered to have different participation levels (e.g., the number of transmitted measurements per unit of time) set by the central planner, which would, in a linear manner, affect their costs (with the coefficients being only known by them). The quality of the service  is a function of the participation levels and the quality of the provided measurements. The central planner relies on the past experiences to empirically track the quality of the measurements provided by each sensor. The sensors report costs and, based on that, the central planner fixes their participation levels and provides payments. The aforementioned schemes differ in several ways from the framework formulated in this paper. Perhaps most importantly, in the framework developed below the central planner does not measure, estimate or directly set the effort expended by the sensors in providing measurements. Rather the central planner seeks to influence the effort invested in the presence of a contract that governs sensor reward. Furthermore, more general cost functions are accommodated and optimality with respect to budget constraints is explored more explicitly.

The framework of this paper also relates to the problem of designing appropriate incentives for agents to truthfully communicate their preferences, most often the parameters of their cost functions, which has a long history in the economics literature on mechanism design; see, e.g.,~\citep{Jackson2003,nisan2007algorithmic}. In those studies, the central planner wishes to make a social decision after discerning from the agents their fixed and privately-held preferences.
This differs from the problem considered here in that the central planner does not seek to recover the private information associated with the effort expended by each sensor, which is not fixed and merely serves as an indicator of measurement quality. In fact, the efforts are functions of the  contracts provided by the central planner to elicit  higher quality estimates. 

In what follows, the outline of the paper is presented. In Section~\ref{sec:problem}, the underlying problem is formulated in game-theoretic terms. Contracts are constructed in Section~\ref{sec:results}. Finally, a numerical example is given in Section~\ref{sec:example} and concluding remarks are presented in Section~\ref{sec:conclusions}. 

\section{Problem Formulation} \label{sec:problem}
Consider the case where a central planner wishes to obtain an accurate measurement of a 
variable\footnote{The scalar nature of the variable of interest is without the loss of generality as entries of a vector can individually be treated in a similar manner.} $x\in\mathbb{R}$ by obtaining measurements from $n$ sensors. The sensors need to expend an effort denoted by $a_i\in\mathbb{R}_{\geq 0}$ to measure the variable $x$. The effort is unknown to the central planner. It determines the quality of the corresponding sensor measurement, which is given by
\begin{align*}
y_i=x+w_i,
\end{align*}
where $w_i$ is a zero-mean random variable with the variance $\mathbb{E}\{w_i^2\}=\eta_i(a_i)$ and $\eta_i:\mathbb{R}_{\geq 0}\rightarrow\mathbb{R}_{\geq 0}$ is an appropriate mapping that captures the return on the effort expended by the sensor. 

\begin{assumption} \label{assum:independent} $\mathbb{E}\{w_iw_j\}=0$ if $i\neq j$.
\end{assumption}

\begin{remark} Note that Assumption~\ref{assum:independent} is satisfied so long as the sensors do not have access to each other's measurements and the estimate constructed by the central planner when reporting their measurements.  This is satisfied, in the setup of this paper, as  the sensors are not allowed to renegotiate their contract with the central planner and update their measurements after the estimate constructed by the central planner is revealed. 
\end{remark}

\begin{assumption} \label{asm:eta} The mapping $\eta_i:\mathbb{R}_{\geq 0}\rightarrow\mathbb{R}_{\geq 0}$ is twice continuously differentiable and strictly decreasing.
\end{assumption}

The pay-off to sensor $i$, for the measurement that it provides to the central planner, is modelled as 
\begin{align*}
C_i(a_i,a_{-i})=\alpha_ip_i-f_i(a_i),
\end{align*} 
where $p_i\in\mathbb{R}_{\geq 0}$ is a potentially stochastic unit of compensation offered by the central planner to sensor $i$ depending on all  reported measurements, $f_i:\mathbb{R}_{\geq 0}\rightarrow\mathbb{R}$ is a mapping that determines the cost to sensor $i$ for investing an effort equal to $a_i$, and $\alpha_i\in\mathbb{R}_{\geq 0}$ is the value-of-compensation\footnote{This can be determined by surveying the sensors or utilizing historical data. Alternatively, the parameter can be eliminated by replacing it with the average value for the society; however, this results in an approximate analysis. } from the perspective of sensor $i$. Also note the game-theoretic notation  $a_{-i}=(a_j)_{j\neq i}$. It is assumed that the sensors deal with the expected cost, that is, they wish to optimize
\begin{align*}
\bar{C}_i(a_i,a_{-i})
&=\mathbb{E}\{C_i(a_i,a_{-i})\}\\
&=\mathbb{E}\{\alpha_ip_i-f_i(a_i)\}\\
&=\alpha_i\mathbb{E}\{p_i\}-f_i(a_i).
\end{align*}
Throughout the paper, rewards take the form $p_i=\pi_i(y_1,\dots,y_n)$ for a given contract $\pi:\mathbb{R}^n\rightarrow \mathbb{R}^n$. The term contract is used because, while the compensation mapping $\pi:\mathbb{R}^n\rightarrow \mathbb{R}^n$ is  agreed prior to expending the effort and gathering data, the level of compensation $p_i$ is determined after the sensors report their measurements. That is, when the sensors agree to participate in this process, they are given a contract for their reports in the future. Moreover, the assumption of dealing with the expectation of the cost is only useful/valid if the sensors provide several reports so that their average returns are well-modelled with the expectation of their cost functions.

\begin{assumption} \label{asm:f} The mapping $f_i:\mathbb{R}_{\geq 0}\rightarrow\mathbb{R}$ is twice continuously differentiable and strictly increasing.
\end{assumption}

\begin{definition}[Contract Game]
A contract game with compensation contract $\pi:\mathbb{R}^n\rightarrow \mathbb{R}^n$ is defined as a tuple $(n,(\mathbb{R}_{\geq 0})_{i=1}^n,(\bar{C}_i)_{i=1}^n)$ that encodes $n$ effort-averse sensors each with the action space $\mathbb{R}_{\geq 0}$ and utility $\bar{C}_i(a_i,a_{-i})=\alpha_i\mathbb{E}\{\pi_i(y_1,\dots,y_n)\}-f_i(a_i)$. The parameters of the game include $(\alpha_i)_{i=1}^n$ as well as the parameters of mappings $\pi$, $(f_i)_{i=1}^n$, and $(\eta_i)_{i=1}^n$.
\end{definition}

\begin{definition}[Contract Equilibrium] \label{def:eq} An action tuple $(a_i^*)_{i=1}^n\in\mathbb{R}_{\geq 0}^n$ constitutes an equilibrium of the contract game if $a_i^*\in\argmax_{a_i\in\mathbb{R}_{\geq 0}} \bar{C}_i(a_i,a_{-i}^*)$ for all $i\in\{1,\dots,n\}$.
\end{definition}

The central planner employs a simple averaging estimator to extract 
\begin{align} \label{eqn:filter}
\hat{x}=\frac{1}{n}(y_1+\dots+y_n).
\end{align}
Note that the central planner cannot implement the least mean square estimator since it cannot directly measure nor set the effort to be invested by the sensors. Asking sensors to communicate the effort that they invest provides scope for malicious behaviors or self-serving decisions, which would ultimately increase the complexity of the central planners.
The quality of the estimate is given by
\begin{align*}
\mathbb{E}\{(x-\hat{x})^2\}
&=\mathbb{E}\bigg\{\bigg(x-\frac{1}{n}(y_1+\dots+y_n)\bigg)^2\bigg\}\\
&=\mathbb{E}\bigg\{\bigg(\frac{1}{n}(w_1+\dots+w_n)\bigg)^2\bigg\}\\
&=\frac{1}{n^2}\sum_{i=1}^n \eta_i(a_i).
\end{align*}
Throughout this paper, it is assumed that the amount of effort $a_i$ expended by sensor $i$ is only known to itself and that this is, implicitly, a function of the devised compensation contract. 

Clearly, having a fixed compensation policy such that $p_i=c$ for all $a_i\in\mathbb{R}_{\geq 0}$ and all $i\in\{1,\dots,n\}$ is not good. This is because, in that case, $\bar{C}_i(a_i,a_{-i})=\alpha_ic-f_i(a_i)$ for all $a_i\in\mathbb{R}_{\geq 0}$. Now, by Assumption~\ref{asm:f}, it can be seen that it is in the best interest of the sensor to select $a_i=0$. This is not the preferable outcome for the central planner in terms of the estimate quality. To fix this issue, it important that $\mathbb{E}\{p_i\}$ becomes a function of $a_i$. This is the topic of the next section.

\section{Compensation Based on Empirical Statistics} \label{sec:results}
Consider the compensation policy
\begin{align}
p_i
&=\pi_i(y_1,\dots,y_n)\nonumber\\
&=\delta_i-\gamma_i(\hat{x}-y_i)^2, \label{eqn:contract}
\end{align}
where $\gamma_i,\delta_i\in\mathbb{R}_{\geq 0}$, $1\leq i\leq n$, are appropriately-selected constants. At first sight, the compensation policy in~\eqref{eqn:contract} might seem restrictive, or  arbitrary. However, at the end of this section, is it established that this policy requires the least total budget for achieving a specified level of performance. Further, for a given budget, it achieves the smallest estimation error variance.
For this compensation contract, by Assumption~\ref{assum:independent}, 
\begin{align*}
\mathbb{E}\{p_i\}
&=\mathbb{E}\{\delta_i-\gamma_i(\hat{x}-y_i)^2\}\\
&=\delta_i-\gamma_i\mathbb{E}\bigg\{\bigg(\frac{1}{n}(y_1+\dots+y_n)-y_i\bigg)^2\bigg\}\\
&=\delta_i-\gamma_i\mathbb{E}\bigg\{\bigg(-\frac{n-1}{n}w_i+\frac{1}{n}\sum_{j\neq i}w_j\bigg)^2\bigg\}\\
&=\delta_i-\gamma_i\bigg(\bigg(\frac{n-1}{n}\bigg)^2\eta_i(a_i)+\frac{1}{n^2}\sum_{j\neq i}\eta_j(a_j)\bigg).
\end{align*}
Consequently
\begin{align} 
\bar{C}_i(a_i,a_{-i})=\alpha_i\delta_i-\bigg[&\alpha_i\gamma_i\bigg(\bigg(\frac{n-1}{n}\bigg)^2\eta_i(a_i)\nonumber\\
&+\frac{1}{n^2}\sum_{j\neq i}\eta_j(a_j)\bigg)+f_i(a_i)\bigg].\label{eqn:expectedcost}
\end{align}
It is is now possible to prove the following intermediate result. 

\begin{proposition} Let
\begin{align}
\mathcal{A}_i=\argmin_{a_i\in\mathbb{R}_{\geq 0}}\bigg[\alpha_i\gamma_i\bigg(\frac{n-1}{n}\bigg)^2\eta_i(a_i)+f_i(a_i)\bigg]. \label{eqn:defmathcalA}
\end{align}
Any action tuple $(a_i^*)_{i=1}^n\in \mathbb{R}_{\geq 0}^n$ is a contract equilibrium if and only if $(a_i^*)_{i=1}^n\in\prod_{i=1}^n \mathcal{A}_i$.
\end{proposition}

\begin{IEEEproof} In view of Definition~\ref{def:eq}, the proof follows immediately by maximizing~\eqref{eqn:expectedcost}. \QED
\end{IEEEproof}

\begin{definition}[Symmetric Contract Game] The contract game is symmetric if $\alpha_i=\alpha$, $\gamma_i=\gamma$, $\delta_i=\delta$, $f_i=f$, and $\eta_i=\eta$ for all $i\in\{1,\dots,n\}$.
\end{definition}

\begin{definition}[Symmetric Contract Equilibrium] A contract equilibrium $(a_i^*)_{i=1}^n$ is symmetric if $a_i^*=a_j^*$ for all $i,j\in\{1,\dots,n\}$.
\end{definition}

\begin{corollary} \label{cor:sym} Any equilibrium of a symmetric contract game is symmetric.
\end{corollary}

\begin{IEEEproof} The proof follows directly from the definition of set $\mathcal{A}_i$ in~\eqref{eqn:defmathcalA} for the case of symmetric contract games. 
\QED \end{IEEEproof}

\begin{remark}[Equilibrium vs. Dominant Strategy] Note that, from the structure of the expected cost function in~\eqref{eqn:expectedcost}, the best action of sensor $i$ is independent of the actions of the other sensors. Therefore, $a_i^*\in\mathcal{A}_i$ is a dominant strategy for the sensor, which is stronger than an equilibrium where the sensors do not deviate given that the others do not deviate as well. This means, even if some sensors are faulty or mistakenly determine their effort, it is in the best interest of each sensor to expend the effort $a_i^*\in\mathcal{A}_i$. This behaviour makes the estimator robust to individual corruptions. Note that this observation however does not exclude the possibility of the sensors being able to improve their compensation by colluding with each other.
\end{remark}

\begin{proposition} \label{prop:existence} A contract equilibrium exists if $\lim_{a_i\rightarrow\infty} f_i(a_i)=\infty$.
\end{proposition}

\begin{IEEEproof} To prove this result, it is established that $\mathcal{A}_i\neq \emptyset$ for all $i$. Let $\xi_i:\mathbb{R}_{\geq 0}\rightarrow \mathbb{R}$ be such that
\begin{align} \label{eqn:defxi}
\xi_i(a_i)=\bigg[\alpha_i\gamma_i\bigg(\frac{n-1}{n}\bigg)^2\eta_i(a_i)+f_i(a_i)\bigg], \,\forall a_i\in\mathbb{R}_{\geq 0}.
\end{align}
It follows that $\lim_{a_i\rightarrow\infty}\xi_i(a_i)\geq \lim_{a_i\rightarrow\infty}f_i(a_i)=+\infty,$ where the inequality follows, by definition, from the property $\eta_i(a_i)\geq 0$ for all $a_i\in\mathbb{R}_{\geq 0}$.
As a result, for all $\Xi\in\mathbb{R}_{\geq 0}$, there exists $A(\Xi)\in\mathbb{R}_{\geq 0}$ such that $\xi_i(a_i)\geq \Xi$ for all $a_i\geq A(\Xi)$. Let $\Xi'\in\mathbb{R}_{\geq 0}$ be an arbitrary real number such that 
$\Xi'>\xi_i(0)=\alpha_i\gamma_i((n-1)^2/n^2)\eta_i(0)+f_i(0).$ Now, define the set $\Omega:=\{a_i\in\mathbb{R}_{\geq 0}\,|\,a_i\leq A(\Xi') \}.$ Clearly, the minimizer of $\xi_i(\cdot)$ belongs to $\Omega$ because $\xi_i(a_i)>\Xi'>\xi_i(0)$ for all $a_i\in\mathbb{R}_{\geq 0}\setminus\Omega$. Hence, $\mathcal{A}_i=\argmin_{a_i\in\Omega}\xi_i(a_i).$ Noting the continuous function $\xi_i(\cdot)$ attains its minimum over the compact set $\Omega$, it follows that $\mathcal{A}_i$ is  non-empty.\QED 
\end{IEEEproof}

Throughout this report, let $f'_i:\mathbb{R}_{\geq 0}\rightarrow\mathbb{R}_{\geq 0}$ and $f''_i:\mathbb{R}_{\geq 0}\rightarrow\mathbb{R}$ denote the first and the second derivatives of $f_i(\cdot)$, respectively; these exists by Assumption~\ref{asm:f}. Similarly, $\eta'_i:\mathbb{R}_{\geq 0}\rightarrow\mathbb{R}_{\leq 0}$ and $\eta''_i:\mathbb{R}_{\geq 0}\rightarrow\mathbb{R}$ denote the first and the second derivatives of $\eta_i(\cdot)$, respectively, which exists by Assumption~\ref{asm:eta}. For some cases, the construction of $\mathcal{A}_i$ can be simplified. The next corollary presents one such case.

\begin{corollary} \label{cor:derivative} Let $\alpha_i\gamma_i[(n-1)^2/n^2]\eta'_i(0)+f'_i(0)<0$ and $\alpha_i\gamma_i[(n-1)^2/n^2]\eta''_i(a)+f''_i(a)\geq 0$ for all $a\in\mathbb{R}_{\geq 0}$. Define
\begin{align*}
\mathcal{A}_i=\bigg\{a_i\in\mathbb{R}_{\geq 0}\,|\, \alpha_i\gamma_i\bigg(\frac{n-1}{n}\bigg)^2\eta'_i(a_i)+f'_i(a_i)=0\bigg\}.
\end{align*}
Any action tuple $(a_i^*)_{i=1}^n\in \mathbb{R}_{\geq 0}^n$ is a contract equilibrium if and only if $(a_i^*)_{i=1}^n\in\prod_{i=1}^n \mathcal{A}_i$.
\end{corollary}

\begin{IEEEproof} For $\xi_i(\cdot)$, defined in~\eqref{eqn:defxi},
$\mathrm{d}^2\xi_i(a_i)/\mathrm{d} a_i^2\geq 0$. Therefore, $\xi_i(\cdot)$ is a convex function. Furthermore,
\begin{align*}
\frac{\mathrm{d}}{\mathrm{d} a_i}\xi_i(a_i)\bigg|_{a_i=0}
&=\alpha_i\gamma_i\bigg(\frac{n-1}{n}\bigg)^2\eta'_i(0)+f'_i(0)<0.
\end{align*}
This shows that $a_i=0$ cannot be a minimizer of $\xi_i(\cdot)$. Therefore, the optimizer belongs to the interior of the set $\mathbb{R}_{\geq 0}$ and, as a result, it should satisfy
\begin{align*}
\frac{\mathrm{d}}{\mathrm{d} a_i}\bigg[\alpha_i\gamma_i\bigg(\frac{n-1}{n}\bigg)^2\eta_i(a_i)+f_i(a_i)\bigg]=0.
\end{align*}
This concludes the proof.
\QED \end{IEEEproof}

\begin{proposition} \label{prop:unique} The contract equilibrium is unique if $\lim_{a\rightarrow\infty} f(a)=\infty$ and $\alpha_i\gamma_i[(n-1)^2/n^2]\eta''_i(a)+f''_i(a)> 0$ for all $a\in\mathbb{R}_{\geq 0}$.
\end{proposition}

\begin{IEEEproof} 
First, note that $\lim_{a\rightarrow\infty} f(a)=+\infty$ guarantees the existence of the contract equilibrium according to Proposition~\ref{prop:existence}. Define the set $\Omega$ and the mapping $\xi_i:\mathbb{R}_{\geq 0}\rightarrow \mathbb{R}$ as in the proof of Proposition~\ref{prop:existence}. 
Note that $\mathrm{d}^2\xi_i(a_i)/\mathrm{d} a_i^2> 0$ and, therefore, $\xi_i(\cdot)$ is a strictly convex function. Therefore, the minimizer of $\xi_i(\cdot)$ over the convex set $\Omega$ is unique.
\QED \end{IEEEproof}

Assume that the conditions of Corollary~\ref{cor:derivative} and Proposition~\ref{prop:unique} hold. Define the mapping $a_i^*:\mathbb{R}_{\geq 0}\rightarrow \mathbb{R}_{\geq 0}$ to be such that
\begin{align} \label{eqn:def:a}
\alpha_i\gamma_i\bigg(\frac{n-1}{n}\bigg)^2\eta'_i(a_i^*(\gamma_i))+f'_i(a_i^*(\gamma_i))=0.
\end{align}
Following the results of Corollary~\ref{cor:derivative} and Proposition~\ref{prop:unique}, the action tuple $(a_i^*(\gamma_i))_{i=1}^n$ is well-defined and is a contract equilibrium for given reward policy parameters $(\gamma_i)_{i=1}^n$.

\begin{proposition} \label{prop:moremoneybetter} Under the assumptions of Corollary~\ref{cor:derivative} and Proposition~\ref{prop:unique}, $a_i^*(\gamma_i)$ is an increasing function of $\gamma_i$.
\end{proposition}

\begin{IEEEproof} First, note that the mapping $a_i^*(\cdot)$ is differentiable; see Theorem~4 in~\citep[p.\,139]{cheney2013analysis}. Taking the derivative of~\eqref{eqn:def:a} with respect to $\gamma_i$ yields
\begin{align*}
0
=&\frac{\mathrm{d}}{\mathrm{d}\gamma_i}\bigg[\alpha_i\gamma_i\bigg(\frac{n-1}{n}\bigg)^2\eta'_i(a_i^*(\gamma_i))+f'_i(a_i^*(\gamma_i)) \bigg]\\
=&\alpha_i\bigg(\frac{n-1}{n}\bigg)^2\eta'_i(a_i^*(\gamma_i))\\
&+\bigg[\alpha_i\gamma_i\bigg(\frac{n-1}{n}\bigg)^2\eta''(a_i^*(\gamma_i))+f''(a_i^*(\gamma_i))\bigg]\frac{\mathrm{d} a_i^*(\gamma_i)}{\mathrm{d}\gamma_i}.
\end{align*}
Consequently
\begin{align*}
\frac{\mathrm{d} a_i^*(\gamma_i)}{\mathrm{d}\gamma_i}=\frac{-\alpha_i[(n-1)^2/n^2]\eta'_i(a_i^*(\gamma_i))}{\alpha_i\gamma_i[(n-1)^2/n^2]\eta''_i(a_i^*(\gamma_i))+f''_i(a_i^*(\gamma_i))}>0,
\end{align*}
which follows from Assumption~\ref{asm:eta} and the inequality $\alpha_i\gamma_i[(n-1)^2/n^2]\eta''_i(a)+f''_i(a)> 0$. 
\QED \end{IEEEproof}

\begin{definition} A compensation policy is \textit{ex ante} individually rational if 
$\bar{C}_i(a_i,a_{-i})\geq 0$, $1\leq i\leq n$, for any contract equilibrium $(a_i)_{i=1}^n\in\prod_{i=1}^n\mathcal{A}_i$.
\end{definition}

Note that, in a liberal society, individually rational compensation policies must be used as, within such a society, the sensors are free to leave the sensing scheme in the next round of contract negotiations or to not participate in the current round of negotiations if there is no hope of receiving compensation in return for efforts. 

\begin{proposition} \label{prop:indvrational} 
The proposed compensation contract in~\eqref{eqn:contract} is \textit{ex ante} individually rational if 
\begin{align*}
\delta_i\geq \gamma_i\bigg(\bigg(\frac{n-1}{n}\bigg)^2\eta_i(a^*_i)+\frac{1}{n^2}\sum_{j\neq i}\eta_j(a^*_j)\bigg)+\frac{1}{\alpha_i}f(a^*_i)
\end{align*}
for all $(a^*_i)_{i=1}^n\in\prod_{i=1}^n\mathcal{A}_i$.
\end{proposition}

\begin{IEEEproof} Note that
\begin{align*}
\bar{C}_i(a^*,a^*_{-i})
&=\alpha_i\delta_i-\bigg[\alpha_i\gamma_i\bigg(\bigg(\frac{n-1}{n}\bigg)^2\eta_i(a^*_i)\\
&\hspace{.8in}+\frac{1}{n^2}\sum_{j\neq i}\eta_j(a^*_j)\bigg)+f_i(a^*_i)\bigg]\\
&=\alpha\bigg[\delta_i-\gamma_i\bigg(\bigg(\frac{n-1}{n}\bigg)^2\eta_i(a^*_i)\\
&\hspace{.8in}+\frac{1}{n^2}\sum_{j\neq i}\eta_j(a^*_j)\bigg)-\frac{1}{\alpha_i}f(a^*_i) \bigg].
\end{align*}
Requiring $\bar{C}_i(a^*,a^*_{-i})\geq 0$ gives the condition in the statement of the proposition. 
\QED \end{IEEEproof}

The condition in Proposition~\ref{prop:indvrational} can be greatly simplified for symmetric contract games.

\begin{corollary} \label{cor:indvrational:symmetric} 
The proposed compensation contract in~\eqref{eqn:contract} is \textit{ex ante} individually rational for a symmetric contract game if 
$\delta\geq \gamma (n-1)\eta(a_i^*)/n+f(a_i^*)/\alpha, \,\forall a_i^*\in\mathcal{A}_i.$
\end{corollary}

\begin{IEEEproof} Following Corollary~\ref{cor:sym}, it is known that the contract equilibrium is symmetric, that is, $a_1^*=\cdots=a_n^*=a^*$. Therefore,
\begin{align*}
\bar{C}_i(a^*,a^*_{-i})
&=\alpha\bigg[\delta-\bigg(\gamma\bigg(\frac{n-1}{n}\bigg)\eta(a^*)+\frac{1}{\alpha}f(a^*)\bigg) \bigg]\geq 0.
\end{align*}
This concludes the proof. 
\QED \end{IEEEproof}

Under the conditions of Corollary~\ref{cor:derivative} and Proposition~\ref{prop:unique}, the total expected budget needed for implementing~\eqref{eqn:contract} at the contract equilibrium $(a_i^*)_{i=1}^n$ when its parameters are set to be individual-rationality is given by
\begin{align}
B
:=&\sum_{i=1}^n\mathbb{E}\{p_i\}\nonumber\\
=&\sum_{i=1}^n \bigg[\delta_i-\gamma_i\bigg(\bigg(\frac{n-1}{n}\bigg)^2\eta_i(a_i^*)\hspace{-.03in}+\hspace{-.03in}\frac{1}{n^2}\sum_{j\neq i}\eta_i(a_j^*)\bigg)\bigg]\label{eqn:budget}\\
\geq &\sum_{i=1}^n f(a_i^*)/\alpha_i,
\label{eqn:budget:lowerbound}
\end{align}
where the inequality in~\eqref{eqn:budget:lowerbound} follows from the condition of Proposition~\ref{prop:indvrational}.
This can be generalized to consideration of all \textit{ex ante} individually rational compensation contracts. 

\begin{proposition} \label{prop:IR} Let $\pi:\mathbb{R}^n\rightarrow\mathbb{R}^n_{\geq 0}$ to be an \textit{ex ante} individually rational compensation contract. Then, the budget for implementing the contract equilibrium $(a_i^*)_{i=1}^n$ is lower bounded as $B\geq \sum_{i=1}^n f_i(a_i^*)/\alpha_i.$
\end{proposition}

\begin{IEEEproof} Individual rationality ensures that $\bar{C}_i(a_i^*,a^*_{-i})\geq 0$ for all $i\in\{1,\dots,n\}$. This implies that $\alpha_i\mathbb{E}\{p_i\}-f_i(a_i^*)\geq 0$ and thus, $B=\sum_{i=1}^n\mathbb{E}\{p_i\}\hspace{-.03in}\geq\hspace{-.03in} \sum_{i=1}^nf_i(a_i^*)/\alpha_i.$ \hfill\QED \end{IEEEproof}

As expected, there is a trade-off between the amount of the budget spent and the quality of the estimate if \textit{ex ante} individually-rational compensation contracts are considered. From this point forward, symmetric contract games are considered. 

\begin{proposition}[Fundamental Budget Requirement]  \label{prop:fund} For any symmetric contract game, the budget for implementing any \textit{ex ante} individually rational compensation contract with estimation quality $\mathbb{E}\{\|x-\hat{x}\|_2^2\}\leq \epsilon$ is lower bounded as $ B\geq nf(\eta^{-1}(\epsilon))/\alpha$.
\end{proposition}

\begin{IEEEproof} By the assumed symmetry and Corollary~\ref{cor:sym}, it follows that $a_1^*=\cdots=a_n^*=a^*$. Notice that
\begin{align*}
\mathbb{E}\{\|x-\hat{x}\|_2^2\}
&=\frac{1}{n}\sum_{i=1}^n \eta(a_i^*)=\eta(a^*).
\end{align*}
Based on Assumption~\ref{asm:eta}, it is known that $\eta^{-1}(\cdot)$ exists. Hence, if  $a^*\geq \eta^{-1}(\epsilon)$ at the contract equilibrium, then the required level of precision is achieved. Following Proposition~\ref{prop:IR}, the budget for implementing this contract equilibrium is lower bounded by
\begin{align*}
B
&\geq \sum_{i=1}^\infty f(a_i^*)/\alpha_i=\frac{n}{\alpha}f(a^*)\geq \frac{n}{\alpha}f(\eta^{-1}(\epsilon)).
\end{align*}
This completes the proof.
\QED \end{IEEEproof}

\begin{corollary}[Fundamental Performance Limit]  \label{prop:performance} For any symmetric contract game, the estimation quality of any \textit{ex ante} individually rational compensation contract with the budget constraint $B\leq \beta$ is lower  bounded as $\mathbb{E}\{\|x-\hat{x}\|_2^2\}\geq \eta(f^{-1}(\beta \alpha/n))$.
\end{corollary}

Note that $\mathbb{E}\{\|x-\hat{x}\|_2^2\}\geq \eta(f^{-1}(\beta \alpha/n))$ is very similar to the Cram\'{e}r-Rao bound for Gaussian estimation problems~\citep{kay1998fundamentals}. This is because of the linear nature of the estimator in~\eqref{eqn:filter}.

Finally, it is shown that the lower bound in Proposition~\ref{prop:fund} is tight and is achieved by the a simple contract of the form~\eqref{eqn:contract}.

\begin{proposition} 
For any symmetric contract game such that $\lim_{a\rightarrow\infty} f(a)=\infty$ and
$\eta''(a)f'(\eta^{-1}(\epsilon))-f''(a)\eta'(\eta^{-1}(\epsilon))>0$ for all $a\in\mathbb{R}_{\geq 0}$, the budget-optimal compensation contract, among the set of all \textit{ex ante} individually-rational compensation contracts guaranteeing a performance level $\mathbb{E}\{\|x-\hat{x}\|_2^2\}\leq \epsilon$, is
\begin{align*}
\pi(y_1,\dots,y_n)=&\bigg[\gamma \frac{n-1}{n}\epsilon+\frac{1}{\alpha}f(\eta^{-1}(\epsilon))\bigg]\\
&-\gamma\bigg(-y_i+\frac{1}{n}\sum_{j=1}^n y_j\bigg)^2,
\end{align*}
with $\gamma=-(n/(n-1))^2f'(\eta^{-1}(\epsilon))/(\alpha \eta'(\eta^{-1}(\epsilon))).$
\end{proposition}

\begin{IEEEproof} First, notice that 
\begin{align*}
\alpha \gamma \bigg(\frac{n-1}{n}\bigg)^2\eta'(\eta^{-1}(\epsilon))+f'(\eta^{-1}(\epsilon))=0,
\end{align*}
which means that $(\eta^{-1}(\epsilon))_{i=1}^n$ is a contract equilibrium. Substituting $\delta=\gamma (n-1)\epsilon/n+f(\eta^{-1}(\epsilon))/\alpha$ and $a_i^*=\eta^{-1}(\epsilon)$, $\forall i$, in~\eqref{eqn:budget} gives
\begin{align*}
B
&=\bigg[\gamma \frac{n-1}{n}\epsilon+\frac{1}{\alpha}f(\eta^{-1}(\epsilon))\bigg] n-\gamma\frac{n-1}{n}\sum_{i=1}^n \eta(\eta^{-1}(\epsilon))\\
&=\frac{n}{\alpha}f(\eta^{-1}(\epsilon)),
\end{align*}
which is the smallest admissible budget according to Proposition~\ref{prop:fund}. 
\QED \end{IEEEproof}

\begin{figure}[t!]
\centering
        \begin{subfigure}[b]{0.9\columnwidth}
\centering
\includegraphics[width=1\columnwidth]{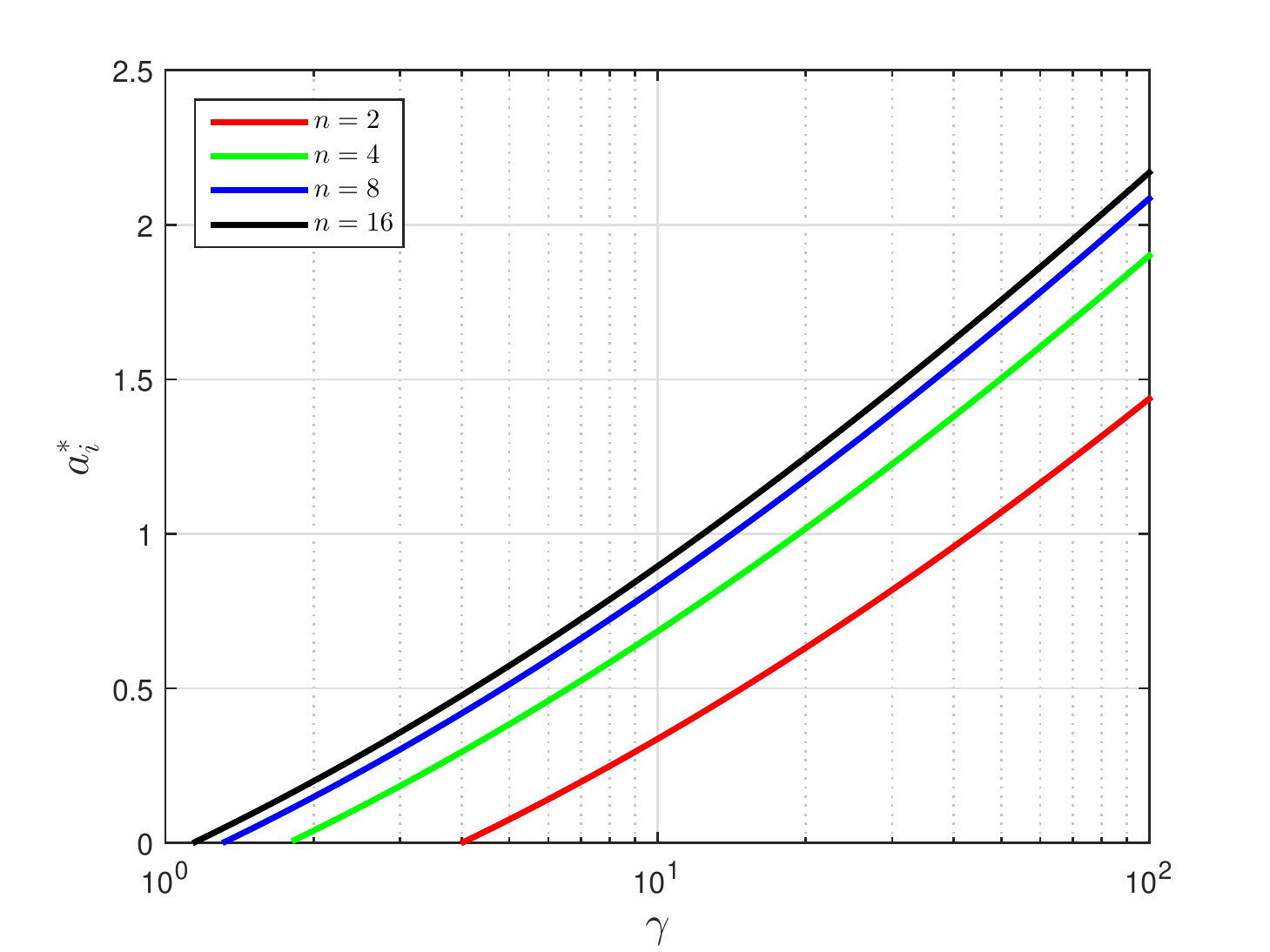}
\caption{\label{fig:1:1}}
        \end{subfigure}
        \begin{subfigure}[b]{0.9\columnwidth}
\centering
\includegraphics[width=1\columnwidth]{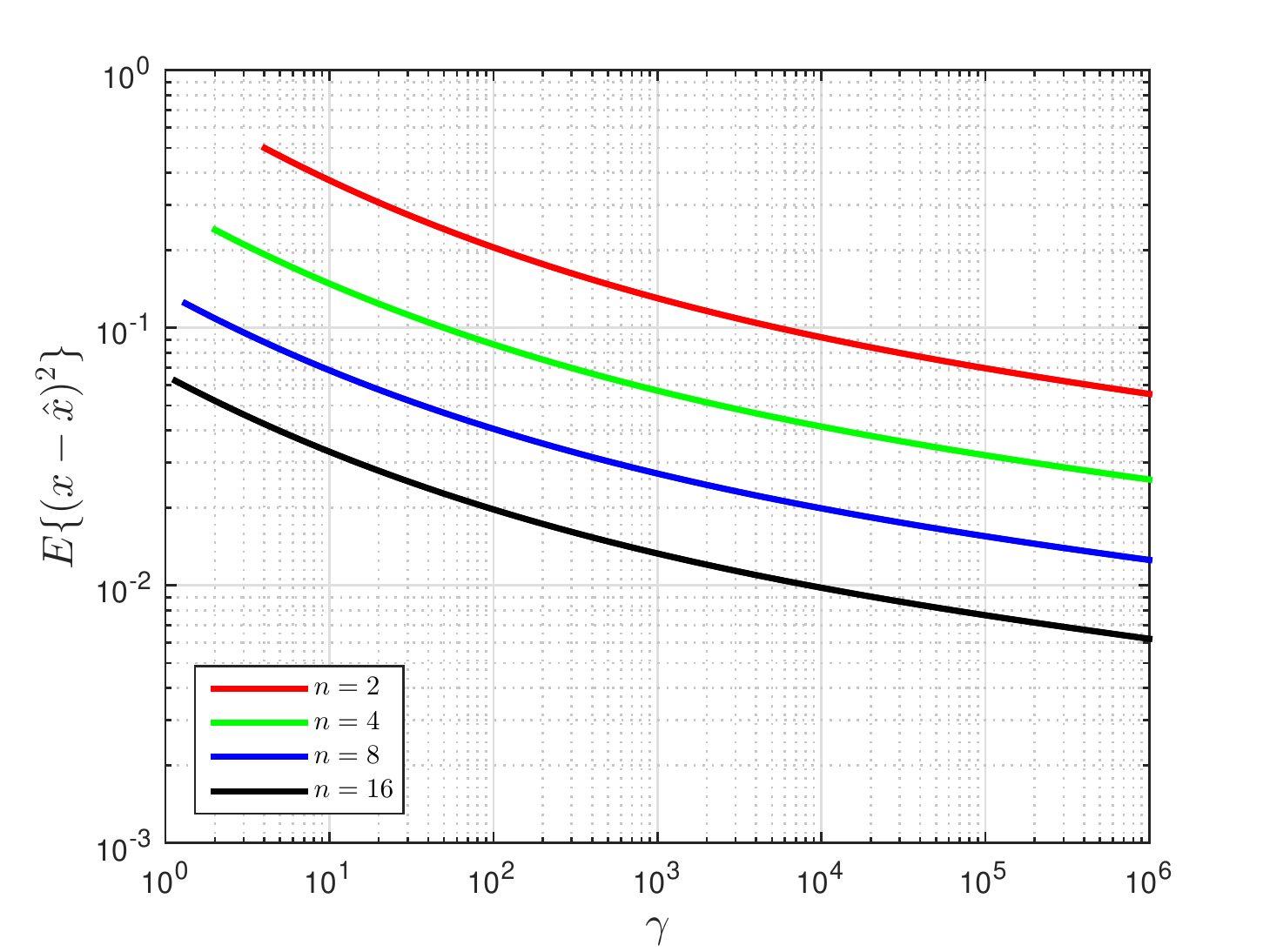}
\caption{\label{fig:1:2}}
        \end{subfigure}
        \begin{subfigure}[b]{0.9\columnwidth}
\centering        
\includegraphics[width=1\columnwidth]{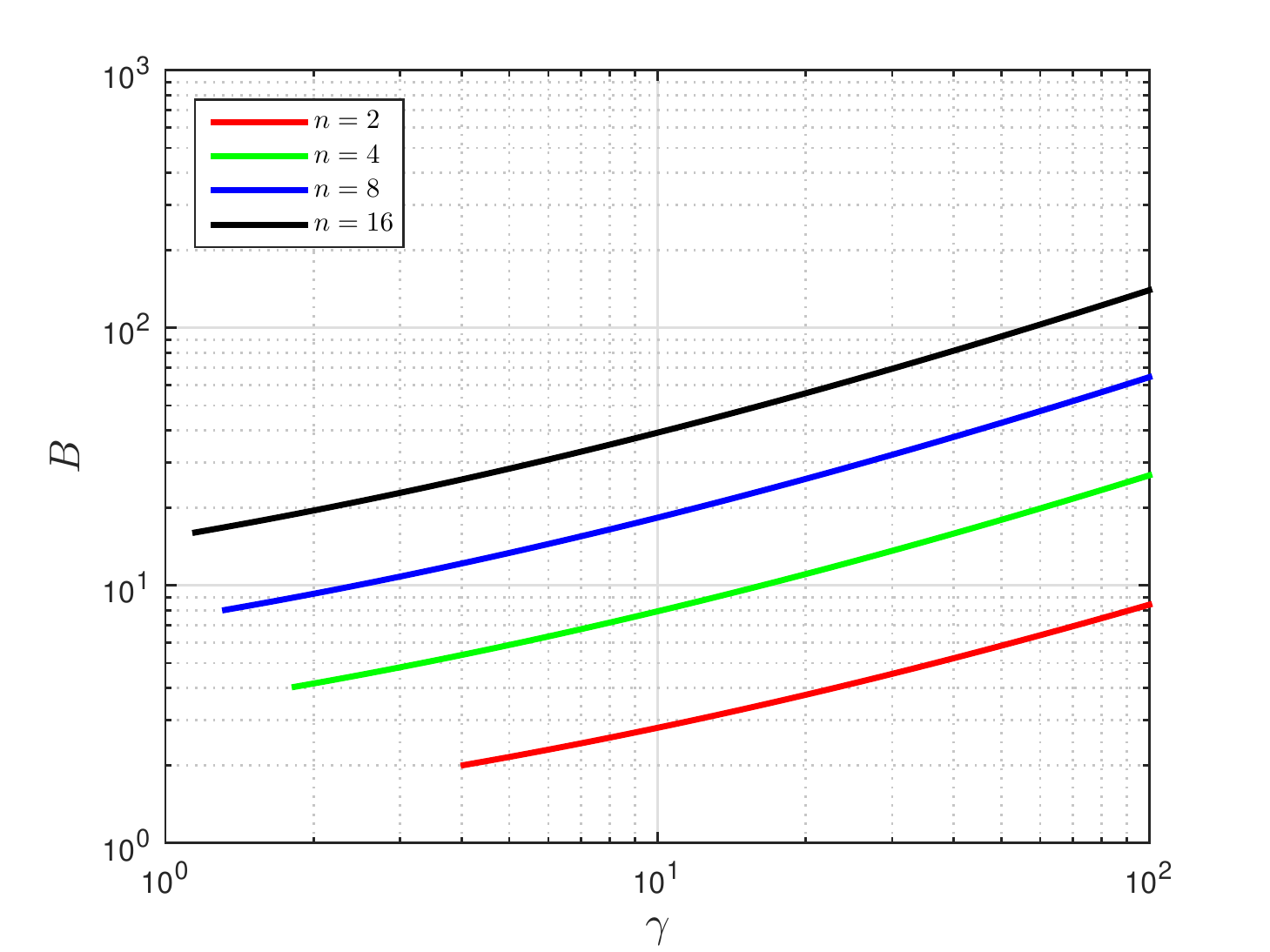}
\caption{\label{fig:1:3}}
        \end{subfigure}
        \caption{\label{fig:1} The effort~(a), the performance~(b), and the budget~(c) at the contract equilibrium for the numerical example as a function of $\gamma$ for various $n$. }
        \vspace{-.1in}
\end{figure}

\section{Numerical Example} \label{sec:example}
Consider a scenario in which sensors are asked to take the measurement of a variable, e.g., the travel time on a road. Noting that their measurements are noisy the sensors may wish to take multiple samples, calculate the average, and transmit that average to the central planner. However, doing this requires the investment of  more time to provide the measurement. Let the effort $a$ for each sensor be proportional to the amount of the time spent refining their measurement. Assume that the cost of that effort is given by $f_i(a)=\exp(\vartheta a)$. This is because, by spending more time, the sensor loses other opportunities for earning money while their estimate does not improve so much as to result in a far superior return from the central planner. Also, let $\eta_i(a)=\varrho/(\varrho+a)$ for a constant $\varrho>0$. This captures the following features: the more time invested, the more samples are gathered and, hence, the more the error in the  measurement can be reduced (which is inversely proportional to the number of the internal samples). Consider a symmetric contract game. Clearly, if $\vartheta-\alpha\gamma[(n-1)^2/n^2]/\varrho< 0$, it follows that
$\alpha\gamma[(n-1)^2/n^2]\eta'_i(0)+f'_i(0)
 < 0.$ Furthermore,
\begin{align*}
\alpha\gamma[(n-1)^2/&n^2]\eta''_i(a)
+f''_i(a) \\
&=2\alpha\gamma[(n-1)^2/n^2]\varrho/(\varrho+a)^3
+
\vartheta^2\exp(\vartheta a)\\
&> 0, \,\forall a\in\mathbb{R}_{\geq 0}.
\end{align*}
Let the compensation policy
$
\pi_i(y_1,\dots,y_n)=\delta-\gamma(\hat{x}-y_i)^2.
$
Upon satisfying $\vartheta-\alpha\gamma[(n-1)^2/n^2]/\varrho< 0$, by Corollary~\ref{cor:derivative} and Proposition~\ref{prop:unique}, the unique contract equilibrium $(a_i^*)_{i=1}^n$ satisfies $-\alpha\gamma((n-1)/n)^2\varrho/(\varrho+a_i^*)^2+\vartheta \exp(\vartheta a_i^*)=0.$
Let us consider the case where $\alpha=1$. Further,  select $\vartheta=1$ and $\varrho=1$. For this contract equilibrium to be individually rational, according to Proposition~\ref{prop:indvrational}, it is sufficient to set $\delta=\gamma (n-1)\eta(a_i^*)/n+f(a_i^*)$.

\begin{figure}
\centering
\includegraphics[width=0.9\columnwidth]{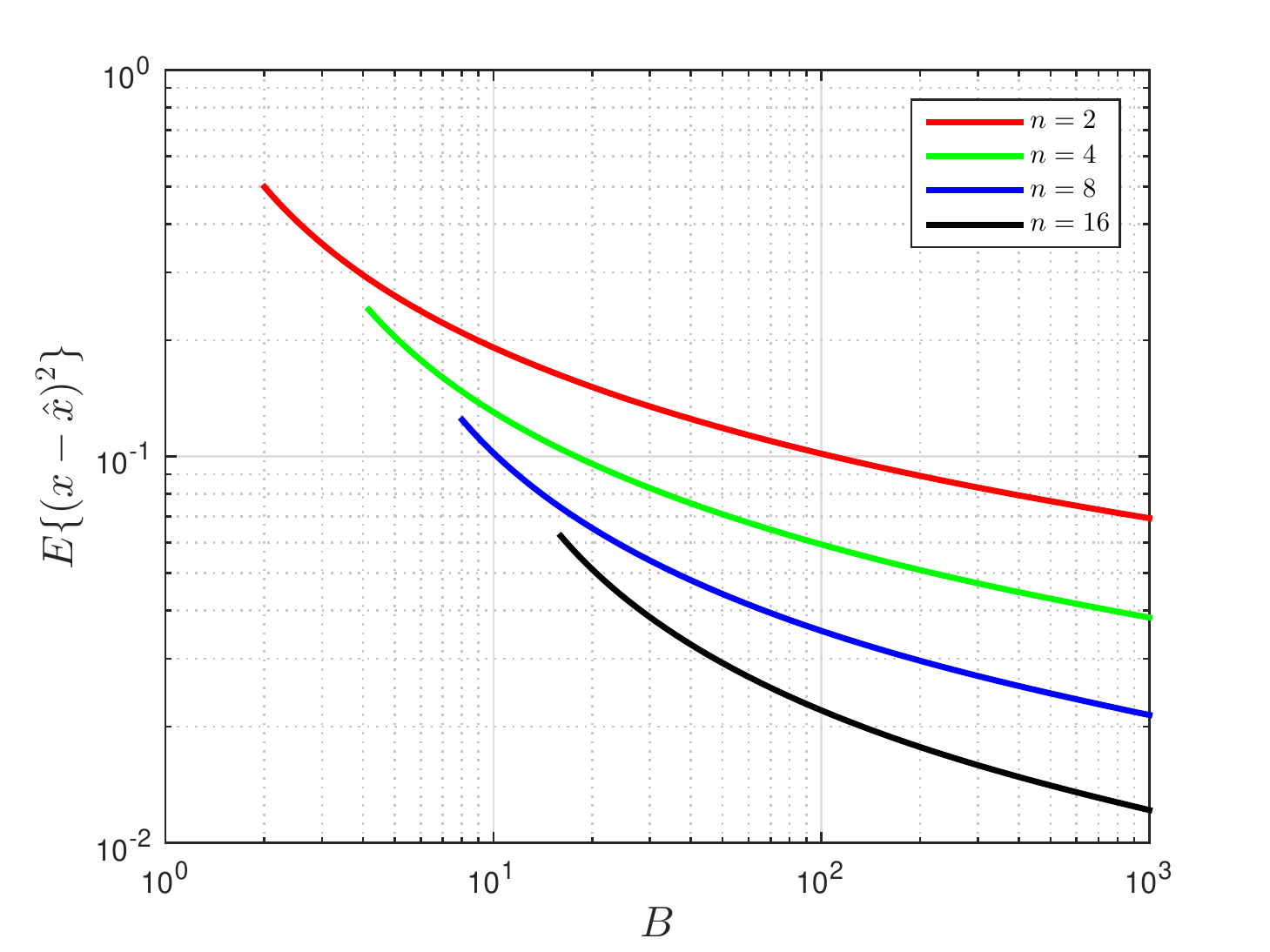}
\caption{\label{fig:2} The performance drawn as a function of the budget for various $n$ at the contract equilibrium for the numerical example.}
\end{figure}

Figure~\ref{fig:1} illustrates the effort expended by the sensors~(a), the performance of the estimator~(b), and the budget required for implementing the contract equilibrium~(c) as a function of $\gamma$ for various numbers of participants $n$. Evidently, the return for adding one extra sensor diminishes with increasing $n$. Figure~\ref{fig:2} shows the required budget versus the quality of the estimate. As expected, with increasing the budget, more accurate measurements are acquired. An interesting observation is that, for a fixed estimation quality, it is better to employ more sensors.

\section{Conclusions and Future Work} \label{sec:conclusions}
The behaviour of effort-averse sensors in response to long term compensations is studied from the perspective of obtaining high-quality measurements. The interaction between the central planner employing an averaging based estimator and sensors via a contract is studied using a game. Conditions for the existence and uniqueness of an equilibrium are identified. Using the characterization of the contract equilibrium, optimal contracts, in terms of the budget, are constructed  for achieving a specified  level of estimation quality.

\bibliographystyle{elsarticle-num}
\bibliography{citation}

\end{document}